\documentclass{article}
\usepackage{amssymb}
\usepackage{amsfonts}

\newtheorem{theorem}{Theorem}

\newtheorem{lemma}[theorem]{Lemma}

\newenvironment{proof}[1][Proof]{\noindent\textbf{#1.} }{\ \rule{0.5em}{0.5em}}
\input{tcilatex}
\begin{document}

\author{Zafer \c{S}iar$^{a}$ and Refik Keskin$^{b}$ \and $\ ^{a}$Bing\"{o}l
University, Mathematics Department, Bing\"{o}l/TURKEY$\ \ $ \and $^{b}$%
Sakarya University, Mathematics Department, Sakarya/TURKEY \and %
zsiar@bingol.edu.tr, rkeskin@sakarya.edu.tr}
\title{On the Diophantine equation $F_{n}-F_{m}=2^{a}$}
\maketitle

\begin{abstract}
In this paper, we solve Diophantine equation in the tittle in nonnegative
integers $m,n,$ and $a$. In order to prove our result, we use lower bounds
for linear forms in logarithms and and a version of the Baker-Davenport
reduction method in diophantine approximation.
\end{abstract}


\section{\protect\bigskip Introduction}

Fibonacci sequence $(F_{n})$ is defined as $F_{0}=0,~F_{1}=1,$ and $%
F_{n}=F_{n-1}+F_{n-2}$ for $n\geq 2.$ The Lucas sequence $(L_{n})$, which is
similar to the Fibonacci sequence, is defined by the same recursive pattern
with initial conditions $L_{0}=2,~L_{1}=1$. The terms of the Fibonacci and
Lucas sequences are called Fibonacci and Lucas numbers, respectively. The
Fibonacci and Lucas numbers for negative indices are defined by $%
F_{-n}=(-1)^{n+1}F_{n}$ and $L_{-n}=(-1)^{n}L_{n}$ \ for $n\geq 1.$ The
Fibonacci and Lucas sequences have many interesting properties and have been
studied in the literature by many researchers. A brief history of Fibonacci
and Lucas sequences one can consult reference \cite{Debnath}. Firstly,
square terms and later perfect powers in the Fibonacci and Lucas sequences
have\ attracted the attention of the researchers. The perfect power in the
Fibonacci and Lucas sequences has been determined in 2006 by Bugeaud,
Mignotte and Siksek in \cite{Bgud} (see Theorem \ref{T2} below). The
Diophantine equation $L_{n}+L_{m}=2^{a}$ $\ $has been tackled in \cite%
{Bravo-1}\ by Bravo and Luca. Two years later, the same authors solved
Diophantine equation $F_{n}+F_{m}=2^{a}$ in \cite{Bravo-2}. Besides, Luca
and Patel, in \cite{Luca1}, found that the Diophantine equation $%
F_{n}-F_{m}=y^{p}$ in integers $(n,m,y,p)$ with $p\geq 2$ has solution
either $\max \left\{ |n|,|m|\right\} \leq 36$ or $y=0$ and $|n|=|m|$ if $%
n\equiv m(\func{mod}2).$ But, it is still an open problem for the case $%
n\not\equiv m(\func{mod}2).$ Motivated by the studies of Bravo and Luca, in
this paper, we consider the Diophantine equation 
\begin{equation}
F_{n}-F_{m}=2^{a}  \label{1.3}
\end{equation}%
in nonnegative integers $m,n,$ and $a.$ We follow the approach and the
method presented in \cite{Bravo-2}. In section 2, we introduce necessary
lemmas and theorems. Then in section 3, we prove our main theorem.

\section{Auxiliary results}

Lately, in many articles, to solve Diophantine equations such as the
equation (\ref{1.3}), authors have used Baker's theory lower bounds for a
nonzero linear form in logarithms of algebraic numbers. Since such bounds
are of crucial importance in effectively solving of Diophantine equations ,
we start with recalling some basic notions from algebraic number theory.

Let $\eta $ be an algebraic number of degree $d$ with minimal polynomial 
\[
a_{0}x^{d}+a_{1}x^{d-1}+...+a_{d}=a_{0}\dprod\limits_{i=1}^{d}\left( X-\eta
^{(i)}\right) \in \mathbb{Z}[x], 
\]%
where the $a_{i}$'s are relatively prime integers with $a_{0}>0$ and $\eta
^{(i)}$'s are conjugates of $\eta .$ Then 
\begin{equation}
h(\eta )=\frac{1}{d}\left( \log a_{0}+\dsum\limits_{i=1}^{d}\log \left( \max
\left\{ |\eta ^{(i)}|,1\right\} \right) \right)  \label{2.1}
\end{equation}%
is called the logarithmic height of $\eta .$ In particularly, if $\eta =a/b$
is a rational number with $\gcd (a,b)=1$ and $b>1,$ then $h(\eta )=\log
\left( \max \left\{ |a|,b\right\} \right) .$

The following properties of logarithmic height are found in many works
stated in references:

\begin{equation}
h(\eta \pm \gamma )\leq h(\eta )+h(\gamma )+\log 2,  \label{2.2}
\end{equation}%
\begin{equation}
h(\eta \gamma ^{\pm 1})\leq h(\eta )+h(\gamma ),  \label{2.3}
\end{equation}%
\begin{equation}
h(\eta ^{s})\leq |s|h(\eta ).  \label{2.4}
\end{equation}%
The following theorem, is deduced from Corollary 2.3 of Matveev \cite{Mtv},
provides a large upper bound for the subscript $n$ in the equation (\ref{1.3}%
) (also see Theorem 9.4 in \cite{Bgud}).

\begin{theorem}
\label{T1} Assume that $\gamma _{1},\gamma _{2},...,\gamma _{t}$ are
positive real algebraic numbers in a real algebraic number field $\mathbb{K}$
of degree $D$, $b_{1},b_{2},...,b_{t}$ are rational integers, and 
\[
\Lambda :=\gamma _{1}^{b_{1}}...\gamma _{t}^{b_{t}}-1 
\]%
is not zero. Then 
\[
|\Lambda |>\exp \left( -1.4\cdot 30^{t+3}\cdot t^{4.5}\cdot D^{2}(1+\log
D)(1+\log B)A_{1}A_{2}...A_{t}\right) , 
\]%
where 
\[
B\geq \max \left\{ |b_{1}|,...,|b_{t}|\right\} , 
\]%
and $A_{i}\geq \max \left\{ Dh(\gamma _{i}),|\log \gamma _{i}|,0.16\right\} $
for all $i=1,...,t.$
\end{theorem}

The following lemma, was proved by Dujella and Peth\H{o} \cite{duj}, is a
variation of a lemma of Baker and Davenport \cite{Baker}. And this lemma
will be used to reduce the upper bound for the subscript $n$ in the equation
(\ref{1.3}). In the following lemma, the function $||\cdot ||$ denotes the
distance from $x$ to the nearest integer, that is, $||x||=\min \left\{
|x-n|:n\in 
%TCIMACRO{\U{2124} }%
%BeginExpansion
\mathbb{Z}
%EndExpansion
\right\} $ for a real number $x.$

\begin{lemma}
\label{L1}Let $M$ be a positive integer, let $p/q$ be a convergent of the
continued fraction of the irrational number $\gamma $ such that $q>6M,$ and
let $A,B,\mu $ be some real numbers with $A>0$ and $B>1.$ Let $\epsilon
:=||\mu q||-M||\gamma q||.$ If $\epsilon >0,$ then there exists no solution
to the inequality 
\[
0<|u\gamma -v+\mu |<AB^{-w}, 
\]%
in positive integers $u,v,$ and $w$ with 
\[
u\leq M\text{ and }w\geq \frac{\log (Aq/\epsilon )}{\log B}. 
\]
\end{lemma}

It is well known that 
\begin{equation}
F_{n}=\dfrac{\alpha ^{n}-\beta ^{n}}{\sqrt{5}}\text{ and }L_{n}=\alpha
^{n}+\beta ^{n},  \label{1.1}
\end{equation}%
where $\alpha =\dfrac{1+\sqrt{5}}{2}$ and $\beta =\dfrac{1-\sqrt{5}}{2},$
which are the roots of the characteristic equations $x^{2}-x-1=0.$ The
relations between Fibonacci and Lucas number, and $\alpha $ are given by 
\begin{equation}
F_{n+1}+F_{n-1}=L_{n}  \label{1.4}
\end{equation}%
and 
\begin{equation}
\alpha ^{n-2}\leq F_{n}\leq \alpha ^{n-1}.  \label{1.2}
\end{equation}%
for $n\geq 1$. The inequality (\ref{1.2}) can be proved by induction. It can
be seen that $1<\alpha <2$ and $-1<\beta <0.$

The following theorem and lemma are given in \cite{Bgud} and \cite{Luca1},
respectively.

\begin{theorem}
\label{T2} The only perfect powers in the Fibonacci sequence are $F_{0}=0,\
F_{1}=F_{2}=1,~F_{6}=8,$ and $F_{12}=144.$ The only perfect powers in the
Lucas sequence are $L_{1}=1$ and $L_{3}=4.$
\end{theorem}

\begin{lemma}
\label{L2}Assume that $n\equiv m(\func{mod}2).$ Then 
\[
F_{n}-F_{m}=\left\{ 
\begin{array}{c}
F_{(n-m)/2}L_{(n+m)/2}\text{ \ \ \ \ \ \ \ if }n\equiv m(\func{mod}4), \\ 
F_{(n+m)/2}L_{(n-m)/2}\text{ \ \ if }n\equiv m+2(\func{mod}4).%
\end{array}%
\right. 
\]
\end{lemma}

\section{Main theorem}

\begin{theorem}
\label{T3}The only solutions of the Diophantine equation (\ref{1.3}) in
nonnegative integers $m<n,$ and $a,$ are given by 
\[
\left( n,m,a\right) \in \left\{ \left( 1,0,0\right) ,\left( 2,0,0\right)
,\left( 3,0,1\right) ,\left( 6,0,3\right) ,\left( 3,1,0\right) ,\left(
4,1,1\right) ,\left( 5,1,2\right) ,\left( 3,2,0\right) \right\} 
\]%
and 
\[
\left( n,m,a\right) \in \left\{ \left( 4,3,0\right) ,\left( 4,2,1\right)
,\left( 5,2,2\right) ,\left( 9,3,5\right) ,\left( 5,4,1\right) ,\left(
7,5,3\right) ,\left( 8,5,4\right) ,\left( 8,7,3\right) \right\} . 
\]
\end{theorem}

%TCIMACRO{\TeXButton{Proof}{\proof}}%
%BeginExpansion
\proof%
%EndExpansion
Assume that the equation (\ref{1.3}) holds. With the help of \textit{%
Mathematica} program, we obtain the solutions in Theorem \ref{T3} for $1\leq
m<n\leq 200.$ This takes a little time. From now on, assume that $n>200$ and 
$n-m\geq 3.$ Now, let us show that $a<n.$ Using (\ref{1.2}), we get the
inequality 
\[
2^{a}=F_{n}-F_{m}<F_{n}<\alpha ^{n}<2^{n}, 
\]%
that is, $a<n.$

On the other hand, rearranging the equation (\ref{1.3}) as $\dfrac{\alpha
^{n}}{\sqrt{5}}-2^{a}=-F_{m}-\dfrac{\beta ^{n}}{\sqrt{5}}$ and taking
absolute values, we obtain 
\[
\left\vert \dfrac{\alpha ^{n}}{\sqrt{5}}-2^{a}\right\vert =\left\vert F_{m}+%
\dfrac{\beta ^{n}}{\sqrt{5}}\right\vert \leq F_{m}+\dfrac{\left\vert \beta
\right\vert ^{n}}{\sqrt{5}}<\alpha ^{m}+\frac{1}{2}
\]%
by (\ref{1.2}). If we divide both sides of the above inequality by $\dfrac{%
\alpha ^{n}}{\sqrt{5}},$ we get 
\begin{equation}
\left\vert 1-2^{a}\alpha ^{-n}\sqrt{5}\right\vert <\frac{4}{\alpha ^{n-m}},
\label{3.1}
\end{equation}%
where we have used the facts that $\alpha ^{-m}<1$ and $n>m.$ Now, let us
apply Theorem \ref{T1} with $\gamma _{1}:=2,~\gamma _{2}:=\alpha ,~\gamma
_{3}:=\sqrt{5}$ and $b_{1}:=a,~b_{2}:=-n,~b_{3}:=1.$ Note that the numbers $%
\gamma _{1},~\gamma _{2},$ and $\gamma _{3}$ are positive real numbers and
elements of the field $\mathbb{K}=Q(\sqrt{5}),$ so $D=2.$ It can be shown
that the number $\Lambda _{1}:=2^{a}\alpha ^{-n}\sqrt{5}-1$ is nonzero. For,
if $\Lambda _{1}=0,$ then we get 
\[
2^{a}=\dfrac{\alpha ^{n}}{\sqrt{5}}=F_{n}+\dfrac{\beta ^{n}}{\sqrt{5}}%
>F_{n}-1>F_{n}-F_{m}=2^{a},
\]%
which is impossible. Moreover, since $h(\gamma _{1})=\log 2=0.6931...,$ $%
h(\gamma _{2})=\dfrac{\log \alpha }{2}=\dfrac{0.4812...}{2},$ and $h(\gamma
_{3})=\log \sqrt{5}=0.8047...$ by (\ref{2.1}), we can take $%
A_{1}:=1.4,~A_{2}:=0.5,$ and $A_{3}:=1.7.$ Also, since $a<n,$ it follows
that $B:=\max \left\{ |a|,|-n|,1\right\} =n.$ Thus, taking into account the
inequality (\ref{3.1}) and using Theorem \ref{T1}, we obtain%
\[
\dfrac{4}{\alpha ^{n-m}}>\left\vert \Lambda _{1}\right\vert >\exp \left(
-1.4\cdot 30^{6}\cdot 3^{4.5}\cdot 2^{2}(1+\log 2)(1+\log n)\left(
1.4\right) \left( 0.5\right) \left( 1.7\right) \right) 
\]%
and so 
\[
(n-m)\log \alpha -\log 4<1.4\cdot 30^{6}\cdot 3^{4.5}\cdot 2^{2}(1+\log
2)(1+\log n)\left( 1.4\right) \left( 0.5\right) \left( 1.7\right) 
\]%
From the last inequality, a quick computation using \textit{Mathematica}
yields to\textit{\ }%
\begin{equation}
(n-m)\log \alpha <2.4\cdot 10^{12}\log n.  \label{3.2}
\end{equation}%
Now, we try to apply Theorem \ref{T1} a second time. Rearranging the
equation (\ref{1.3}) as $\dfrac{\alpha ^{n}}{\sqrt{5}}-\dfrac{\alpha ^{m}}{%
\sqrt{5}}-2^{a}=\dfrac{\beta ^{n}}{\sqrt{5}}-\dfrac{\beta ^{m}}{\sqrt{5}}$
and taking absolute values in here, we obtain 
\[
\left\vert \dfrac{\alpha ^{n}(1-\alpha ^{m-n})}{\sqrt{5}}-2^{a}\right\vert
=\left\vert \dfrac{\beta ^{n}}{\sqrt{5}}-\dfrac{\beta ^{m}}{\sqrt{5}}%
\right\vert \leq \dfrac{\left\vert \beta \right\vert ^{n}+\left\vert \beta
\right\vert ^{m}}{\sqrt{5}}<\frac{1}{3},
\]%
where we used the fact that $\left\vert \beta \right\vert ^{n}+\left\vert
\beta \right\vert ^{m}<2/3$ for $n>200.$ Dividing both sides of the above
inequality by $\dfrac{\alpha ^{n}(1-\alpha ^{n-m})}{\sqrt{5}},$ we get 
\begin{equation}
\left\vert 1-2^{a}\alpha ^{-n}\sqrt{5}(1-\alpha ^{m-n})^{-1}\right\vert <%
\frac{\sqrt{5}\alpha ^{-n}(1-\alpha ^{m-n})^{-1}}{3}.  \label{3.3}
\end{equation}%
Since%
\[
\alpha ^{m-n}=\frac{1}{\alpha ^{n-m}}<\frac{1}{\alpha }<\frac{2}{3},
\]%
it is seen that%
\[
1-\alpha ^{m-n}>1-\frac{2}{3}=\frac{1}{3},
\]%
and therefore%
\[
\frac{1}{1-\alpha ^{m-n}}<3.
\]%
Then from (\ref{3.3}), it follows that 
\begin{equation}
\left\vert 1-2^{a}\alpha ^{-n}\sqrt{5}(1-\alpha ^{m-n})^{-1}\right\vert <%
\frac{3}{\alpha ^{n}}.  \label{3.a}
\end{equation}%
Thus, taking $\gamma _{1}:=2,~\gamma _{2}:=\alpha ,~\gamma _{3}:=\sqrt{5}%
(1-\alpha ^{m-n})^{-1}$ and $b_{1}:=a,~b_{2}:=-n,~b_{3}:=1,$ we can apply
Theorem \ref{T1}. As one can see that, the numbers $\gamma _{1},~\gamma _{2},
$ and $\gamma _{3}$ are positive real numbers and elements of the field $%
\mathbb{K}=Q(\sqrt{5}),$ so $D=2.$ Since%
\[
\dfrac{\alpha ^{n}}{\sqrt{5}}-\dfrac{\alpha ^{m}}{\sqrt{5}}=F_{n}+\dfrac{%
\beta ^{n}}{\sqrt{5}}-F_{m}-\dfrac{\beta ^{m}}{\sqrt{5}}\neq 2^{a}
\]%
for $n>m,$ the number $\Lambda _{2}:=2^{a}\alpha ^{-n}\sqrt{5}(1-\alpha
^{m-n})^{-1}-1$ is nonzero. Similarly, since $h(\gamma _{1})=\log
2=0.6931...,$ and $h(\gamma _{2})=\dfrac{\log \alpha }{2}=\dfrac{0.4812...}{2%
}$ by (\ref{2.1}), we can take $A_{1}:=1.4$ and $~A_{2}:=0.5.$ Besides,
using (\ref{2.2}), (\ref{2.3}), and (\ref{2.4}), we get that $h(\gamma
_{3})\leq \log 2\sqrt{5}+(n-m)\dfrac{\log \alpha }{2},$ and so we can take $%
A_{3}:=\log 20+(n-m)\log \alpha .$ Also, since $a<n,$ it follows that $%
B:=\max \left\{ |a|,|-n|,1\right\} =n.$ Thus, taking into account the
inequality (\ref{3.a}) and using Theorem \ref{T1}, we obtain%
\[
\dfrac{3}{\alpha ^{n}}>\left\vert \Lambda _{2}\right\vert >\exp (-C)(1+\log
2)(1+\log n)\left( 1.4\right) \left( 0.5\right) \left( \log 20+(n-m)\log
\alpha \right) 
\]%
or 
\[
n\log \alpha -\log 3<C(1+\log 2)(1+\log n)\left( 1.4\right) \left(
0.5\right) \left( \log 20+(n-m)\log \alpha \right) ,
\]%
where $C=1.4\cdot 30^{6}\cdot 3^{4.5}\cdot 2^{2}.$ Inserting the inequality (%
\ref{3.3}) into the last inequality, we get 
\begin{equation}
n\log \alpha -\log 3<C(1+\log 2)(1+\log n)\left( 1.4\right) \left(
0.5\right) \left( \log 20+2.4\cdot 10^{12}\log n\right)   \label{3.4}
\end{equation}%
and so $n<2.91\cdot 10^{28}.$

Now, let us try to reduce the upper bound on $n$ applying Lemma \ref{L1} two
times. Let 
\[
z_{1}:=a\log 2-n\log \alpha +\log \sqrt{5}. 
\]%
Then 
\[
\left\vert 1-e^{z_{1}}\right\vert <\frac{4}{\alpha ^{n-m}} 
\]%
by (\ref{3.1}). The inequality 
\[
\dfrac{\alpha ^{n}}{\sqrt{5}}=F_{n}+\dfrac{\beta ^{n}}{\sqrt{5}}%
>F_{n}-1>F_{n}-F_{m}=2^{a} 
\]%
implies that $z_{1}<0.$ In that case, since $\dfrac{4}{\alpha ^{n-m}}<0.95$
for $n-m\geq 3,$ it follows that $e^{|z_{1}|}<20.$ Hence, we get 
\[
0<|z_{1}|<e^{|z_{1}|}-1\leq e^{|z_{1}|}\left\vert 1-e^{z_{1}}\right\vert <%
\dfrac{80}{\alpha ^{n-m}}, 
\]%
or 
\[
0<|a\log 2-n\log \alpha +\log \sqrt{5}|<\dfrac{80}{\alpha ^{n-m}}. 
\]%
Dividing this inequality by $\log \alpha ,$ we get 
\begin{equation}
0<|a\left( \frac{\log 2}{\log \alpha }\right) -n+\frac{\log \sqrt{5}}{\log
\alpha }|<50\cdot \alpha ^{-(n-m)}.  \label{3.5}
\end{equation}%
Now, we can apply Lemma \ref{L1}. Put 
\[
\gamma :=\dfrac{\log 2}{\log \alpha }\notin 
%TCIMACRO{\U{211a} }%
%BeginExpansion
\mathbb{Q}
%EndExpansion
,~\mu :=\dfrac{\log \sqrt{5}}{\log \alpha },~A:=50,~B:=\alpha \text{, and }%
w:=n-m. 
\]%
Taking $M:=2.91\cdot 10^{28},$ we found that $q_{64},$ the denominator of
the $64^{th}$ convergent of $\gamma $ exceeds $6M.$ Furthermore, 
\[
\epsilon =||\mu q_{64}||-M||\gamma q_{64}||\geq 0.184. 
\]%
Thus, we can say that the inequality (\ref{3.5}) has no solutions for%
\[
n-m\geq \dfrac{\log \left( Aq_{64}/\epsilon \right) }{\log B}. 
\]%
A computer search with \textit{Mathematica} yields to $n-m\geq 146.408.$ So $%
n-m\leq 146.$ Substituting this upper bound for $n-m$ into (\ref{3.4}), we
obtain $n<7.56\cdot 10^{15}.$

Now, let us apply again Lemma \ref{L1} to reduce a little bit the upper
bound on $n.$ Let 
\[
z_{2}:=a\log 2-n\log \alpha +\log \left( \sqrt{5}(1-\alpha
^{m-n})^{-1}\right) .
\]%
In this case, 
\[
\left\vert 1-e^{z_{2}}\right\vert <\frac{3}{\alpha ^{n}}
\]%
by (\ref{3.3}). It is seen that $\dfrac{3}{\alpha ^{n}}<\dfrac{1}{2}.$ If $%
z_{2}>0,$ then $0<z_{2}<e^{z_{2}}-1<\dfrac{3}{\alpha ^{n}}.$ If $z_{2}<0,$
then $\left\vert 1-e^{z_{2}}\right\vert =1-e^{z_{2}}<\dfrac{3}{\alpha ^{n}}<%
\dfrac{1}{2}.$ From this, we get $e^{\left\vert z_{2}\right\vert }<2$ and
therefore 
\[
0<\left\vert z_{2}\right\vert <e^{\left\vert z_{2}\right\vert }-1\leq
e^{\left\vert z_{2}\right\vert }\left\vert 1-e^{z_{2}}\right\vert <\dfrac{6}{%
\alpha ^{n}}.
\]%
In any case, the inequality 
\[
0<\left\vert z_{2}\right\vert <\dfrac{6}{\alpha ^{n}}
\]%
is true. That is, 
\[
0<\left\vert a\log 2-n\log \alpha +\log \left( \sqrt{5}(1-\alpha
^{m-n})^{-1}\right) \right\vert <\dfrac{6}{\alpha ^{n}}.
\]%
Dividing both sides of the above inequality by $\log \alpha ,$ we get 
\begin{equation}
0<\left\vert a\left( \frac{\log 2}{\log \alpha }\right) -n+\frac{\log \left( 
\sqrt{5}(1-\alpha ^{m-n})^{-1}\right) }{\log \alpha }\right\vert <13\cdot
\alpha ^{-n}.  \label{3.6}
\end{equation}%
Putting $\gamma :=\dfrac{\log 2}{\log \alpha }$ and taking $M:=7.6\cdot
10^{15},$ we found that $q_{44},$ the denominator of the $44^{th}$
convergent of $\gamma $ exceeds $6M.$ Also, taking%
\[
\mu :=\dfrac{\log \left( \sqrt{5}(1-\alpha ^{m-n})^{-1}\right) }{\log \alpha 
}
\]%
with $n-m\in \left[ 3,146\right] $ except for $n-m\neq 4,12,$ a quick
computation using \textit{Mathematica} gives us the inequality 
\[
\epsilon =||\mu q_{44}||-M||\gamma q_{44}||\geq 0.49939.
\]%
Let $A:=13,~B:=\alpha ,$ and $w:=n$ in Lemma \ref{L1}. Thus, with the help
of \textit{Mathematica}, we can say that the inequality (\ref{3.6}) has no
solution for $n\geq 98.1915$ with $n-m\neq 4,12.$ In that case $n\leq 98.$
This contradicts our assumption that $n>200.$ Thus, we have to consider the
cases $n-m=1,2,4,$ and $12$ to complete the proof. \ If $n-m=1,$ then we
have the equation $2^{a}=$ $F_{m+1}-F_{m}=F_{m-1},$ which implies that $%
\left( n,m,a\right) =\left( 3,2,0\right) ,\left( 4,3,0\right) ,\left(
8,7,3\right) .$ If $n-m=2,$ then we have the equation $%
2^{a}=F_{m+2}-F_{m}=F_{m+1},$ which implies that $\left( n,m,a\right)
=\left( 3,1,0\right) .$ If $n-m=4,$ then we have the equation $%
2^{a}=F_{m+4}-F_{m}=F_{m+3}+F_{m+1}=L_{m+2}$ by (\ref{1.4}). By Theorem \ref%
{T2}, this is only possible for $m=1,$ which implies that $n=5$ and $a=2.$
Now, assume that $n-m=12.$ Then, we have the equation $F_{m+12}-F_{m}=2^{a}.$
since $m+12\equiv m(\func{mod}4),$ it follows that $%
2^{a}=F_{m+12}-F_{m}=F_{6}L_{m+6}$\ by Lemma \ref{L2}. This implies that $%
L_{m+6}=2^{a-3},$ which is impossible by Theorem \ref{T2} since $m>0.$

\end{document}